\documentclass{amsart}
\usepackage{setspace}

\usepackage{amsmath}
\usepackage{amsthm}
\usepackage{amssymb}
\usepackage{enumerate}

\newtheorem{theorem}{Theorem}[section]
\newtheorem{lemma}[theorem]{Lemma}
\newtheorem{proposition}[theorem]{Proposition}
\newtheorem{corollary}[theorem]{Corollary}
\newtheorem{remark}[theorem]{Remark}

\numberwithin{equation}{section}

\begin{document}

\title[Division quaternion algebras over some cyclotomic fields]
{Division quaternion algebras over some cyclotomic fields}

\author{Diana Savin}
\address{Faculty of Mathematics and Computer Science, Transilvania University\\
 Iuliu Maniu street 50, Bra\c{s}ov 500091, Romania\\  https://www.unitbv.ro/en/contact/search-in-the-unitbv-community/4292-diana-savin.html}
\email{diana.savin@unitbv.ro; dianet72@yahoo.com}

\date{}

\begin{abstract}
 Let $p_{1}, p_{2}$ be two distinct prime integers, let  $n$ be a positive integer, $n$$\geq 3$ and let $\xi_{n} $ be a primitive
root of order $n$ of the unity. In this paper we obtain a complete characterization  for a quaternion algebra $H\left(p_{1}, p_{2}\right)$ to be a division algebra over
the $n$th cyclotomic field $\mathbb{Q}\left(\xi_{n}\right)$, when $n$$\in$$\left\{3,4,6,7,8,9,11,12\right\}$ and also we obtain a characterization  for a quaternion algebra $H\left(p_{1}, p_{2}\right)$ to be a division algebra over
the $n$th cyclotomic field $\mathbb{Q}\left(\xi_{n}\right)$, when $n$$\in$$\left\{5,10\right\}$. In the 4th section we obtain a complete characterization for a quaternion algebra $H_{\mathbb{Q}\left(\xi_{n}\right)}\left(p_{1}, p_{2}\right)$ 
to be a division algebra, when $n=l^{k},$ with $l$ a prime integer, $l\equiv 3$ (mod $4$) and $k$ a positive integer. In the last section  of this article we obtain a complete characterization for a quaternion algebra $H_{\mathbb{Q}\left(\xi_{l}\right)}\left(p_{1}, p_{2}\right)$ 
to be a division algebra, when $l$ is a Fermat prime number.
\end{abstract}

\maketitle

\section{Introduction}
\begin{center}
\end{center}

Let $L$ be a field with $char\left(L\right)\neq2$ and let $a, b$$\in$$L\backslash \{0\}.$ \textit{The generalized quaternion algebra over the field} $L$ is an associative algebra $H_{L}\left(\alpha, \beta\right)$ with a basis $\{1,e_{1},e_{2},e_{3}\}$  and the following multiplication:
\begin{equation*}
\begin{tabular}{c||c|c|c|c|}
$\cdot $ & $1$ & $e_{1}$ & $e_{2}$ & $e_{3}$ \\ \hline\hline
$1$ & $1$ & $e_{1}$ & $e_{2}$ & $e_{3}$ \\ \hline
$e_{1}$ & $e_{1}$ & $\alpha$ & $e_{3}$ & $\alpha e_{3}$ \\ \hline
$e_{2}$ & $e_{2}$ & $-e_{3}$ & $\beta$ & $-\beta e_{1}$ \\ \hline
$e_{3}$ & $e_{3}$ & $-\alpha e_{2}$ & $\beta e_{1}$ & $-\alpha \beta$ \\ \hline
\end{tabular}%
\end{equation*}%
Let $x$$\in$$H_{L}\left(\alpha, \beta\right).$ $x=x_{0}\cdot 1+x_{1}e_{1}+x_{2}e_{2}+x_{3}e_{3},$
with $x_{i}\in L,$ $\left(\forall\right)$$i\in \{0,1,2,3\}$ and $\overline{x}$ is the conjugate of $x,$ $\overline{x}=x_{0}\cdot 1-x_{1}e_{1}-x_{2}e_{2}-x_{3}e_{3},$ the norm of $x$ is $\boldsymbol{n}\left(
x\right) =x\cdot \overline{x}=x_{0}^{2}- \alpha x_{1}^{2}- \beta x_{2}^{2}+ \alpha \beta x_{3}^{2}.$ Quaternion algebras are non-commutative algebras and also central simple algebras over the field $L,$  of dimension $4$ over $L.$ Moreover, quaternion algebras are Clifford algebras.\\
The results about the quaternion algebras have many connections with algebraic number theory (especially with the ramification theory in algebraic number fields) and many applications in combinatorics and signal theory. \\
The ramification theory in algebraic number fields, described so beautifully for cyclotomic fields in Washington's book (\cite{Washington}) and for various algebraic number fields in Ribenboim's book (\cite{ribenboim2}),  Ireland and Rosen's book (\cite{ireland}) is used in the research of many  subjects in the number theory or the associative algebras (see \cite{Min_Sav}, ...).\\
Different criterias are known for a quaternion algebra to split (see \cite{lam}, \cite{Voight}, \cite{Voight1}, \cite{alsina}).
In the specialized literature, several articles have been written about conditions in which various quaternion algebras are with division or split over various  algebraic number fields. For example, in \cite{alsina}, explicit conditions were found for quaternion algebras of type $H\left(p_{1}, p_{2}\right)$ (where $p_{1}, p_{2}$ are prime integers) to be division algebras or to be split over the field of rationals numbers.\\
In \cite{savin2016} and \cite{savin2017} we studied some division or split quaternion algebras over the quadratic field  $\mathbb{Q}\left(i\right).$  After this we determined a sufficient condition for a quaternion algebra $H\left(p,q\right)$ (where $p$ and $q$ are prime positive integers) to split over a quadratic field  $\mathbb{Q}\left(\sqrt {d}\right).$ In \cite{astz} we found explicit necessary and sufficient conditions for a quaternion algebra $H\left(p,q\right)$ (where $p$ and $q$ are prime positive integers) to split or to be with division over a quadratic field  $\mathbb{Q}\left(\sqrt {d}\right).$ In \cite{astz2} we found explicit necessary and sufficient conditions for a quaternion algebra  to be a division algebra over a composite  quadratic number fields. In \cite{astz1} we found a complete
characterization of division quaternion algebras $H\left(p,q\right)$  over $K,$ where $K$ is either a dihedral extension of $\mathbb{Q}$
 of prime degree $l$ over a quadratic field $F$, or an abelian $l$
extension of $F$ unramified over $F$ whenever $l$ divides the class number of $F.$\\
In this paper we study quaternion algebras $H\left(p_{1}, p_{2}\right)$ (where $p_{1}$ and $p_{2}$ are prime integers) over a cyclotomic field $\mathbb{Q}\left(\xi_{n}\right),$ where $n$$\in$$\left\{3,4,5,6,7,8,9,10,11,12\right\}$ (in section 3), respectively $n=l^{k},$ with $l$ a prime integer,
$l\equiv 3$ (mod $4$) and $k$ a positive integer. Finally we study quaternion algebras $H\left(p_{1}, p_{2}\right)$ (where  $p_{1}$ and $p_{2}$ are prime integers) over the Kummer field $\mathbb{Q}\left(\xi_{l^{k}}, \sqrt[l^{k}]{\alpha}\right),$ where $n=l^{k},$ with $l$ a prime integer,
$l\equiv 3$ (mod $4$) and $k$ a positive integer (in section 4).
\begin{equation*}
\end{equation*}

\section{ Preliminaries}
\label{preliminari}

We recall  now some results about quadratic fields, cyclotomic fields and quaternions algebras which we shall use for to prove our results.

  In the following theorem it is presented the decomposition of a prime in the ring of algebraic integers of a quadratic filed.
\begin{theorem}
\label{theorem2.1}
 (\cite{ireland}).
 \textit{Let} $d\neq 0,1$ \textit{be a free squares integer.} \textit{Let} $\mathcal{O}_{K}$ \textit{be the ring of integers of the quadratic field} $K=\mathbb{Q}\left(\sqrt{d}\right)$ \textit{and} $\Delta_{K}$ \textit{be the discriminant of} $K.$ \textit{Let} $p$ \textit{be an odd prime integer. Then, we have}:\\
i)  $p$ \textit{is ramified in} $\mathcal{O}_{K},$ \textit{if and only if} $p|\Delta_{K},$  \textit{exactly} $p\mathcal{O}_{K}=\left(p,\sqrt{d}\right)^{2};$\\
ii)  $p$ \textit{splits completely  in} $\mathcal{O}_{K}$ \textit{if and only if}  \textit{the Legendre's symbol} $\left(\frac{\Delta_{K}}{p}\right)=1,$  \textit{exactly} $p\mathcal{O}_{K}=P\cdot \overline{P},$ \textit{where} $P, \overline{P}$ \textit{are prime ideals in} $\mathcal{O}_{K},$ $P\neq \overline{P};$\\
iii)  $p$ \textit{is inert in} $\mathcal{O}_{K}$ \textit{if and only if}   \textit{the Legendre's symbol} $\left(\frac{\Delta_{K}}{p}\right)=-1$  \\
iv) $d\equiv 2$ (\textit{mod} $4$) \textit{if and only if} $2\mathcal{O}_{K}=\left(2,\sqrt{d}\right)^{2};$\\
v) $d\equiv 3$ (\textit{mod} $4$) \textit{if and only if} $2\mathcal{O}_{K}=\left(2,1+\sqrt{d}\right)^{2};$\\
vi) $d\equiv 1$ (\textit{mod} $8$) \textit{if and only if} $2$ \textit{is totally split in} $\mathcal{O}_{K},$ \textit{exactly} $2\mathcal{O}_{K}=P\cdot P^{'},$ \textit{with} $P, P^{'}$
\textit{are prime ideals in} $\mathcal{O}_{K},$ $P=\left(2, \frac{1+\sqrt{d}}{2}\right)$ $P\neq P^{'};$\\
vii) $d\equiv 5$ (\textit{mod} $8$) \textit{if and only if} $2\mathcal{O}_{K}$ \textit{is a prime ideal in} $\mathcal{O}_{K}.$
\end{theorem}

We recall some of the properties of Legendre's symbol, which we will use in the proof of our results.
\begin{proposition}
\label{twodotfthree}
(\cite{ribenboim2}). Let $p$ and $q$ be two distinct odd prime positive integers and let $\left(\frac{\cdot}{p}\right)$ be the Legendre's symbol. Then, the following statements are true:
 \begin{enumerate}[\rm1.]
 \item
$\left(\frac{-1}{p}\right)=\left(-1\right)^{\frac{p-1}{2}}$;
 \item
$\left(\frac{2}{p}\right)=\left(-1\right)^{\frac{p^{2}-1}{8}}$;
\item
If $a,b$ are two integers, $a\equiv b$ (mod $p$), then $\left(\frac{a}{p}\right)=\left(\frac{b}{p}\right)$;
 \item
The Quadratic reciprocity law:
$\left(\frac{p}{q}\right)\cdot \left(\frac{q}{p}\right)=\left(-1\right)^{\frac{p-1}{2}\cdot \frac{q-1}{2}}$.

 \end{enumerate}
\end{proposition}

 Let $\xi_{n}$ be a primitive root of order $n$ of the unity, where $n$ is an integer, $n\geq 3$ and let $\mathbb{Q}\left(\xi_{n}\right)$ be the $n$th cyclotomic field. The field extension $\mathbb{Q}\subset \mathbb{Q}\left(\xi_{n}\right)$ is a Galois extension of degree $\varphi\left(n\right)$ (see \cite{ribenboim2},  \cite{Stefanescu}). Also, the ring of algebraic integers of the cyclotomic field $\mathbb{Q}\left(\xi_{n}\right)$ is $\mathbb{Z}\left[\xi_{n}\right]$ (see \cite{ribenboim2}).\\
The decomposition of a prime in the ring of algebraic integers of a cyclotomic field is presented in the following theorem (\cite{ireland}, \cite{Washington}).

\begin{theorem}
\label{theorem2.2}
 \textit{Let }$n$\textit{\ be a positive integer, }%
$n\geq 3$\textit{\ and let} $\xi_{n}$\textit{\ be a primitive root of order} $n$ of the unity.
 Let  $\mathbb{Q}\left(\xi_{n}\right)$ be the $n$th cyclotomic field. If $p$\textit{\ is a prime
positive integer,} $n$\textit{\ is not divisible with }$p$\textit{\ and }$f$%
\textit{\ is the smallest positive integer such that} $p^{f^{'}}\equiv 1$  (mod $n$),
\textit{then we have }$p\mathbb{Z}\left[\xi_{n}\right]=P_{1}P_{2}....P_{g},$ 
\textit{where }$g=\frac{\varphi \left( n\right) }{f^{'}},\varphi $\textit{\ is
the Euler's function and }$P_{j},\,j=1,...,g$\textit{\ are different prime
ideals in the ring }$\mathbb{Z}[\xi ].$
\end{theorem}

\begin{corollary}
\label{twodotnine}
(\cite{Washington}) 
\textit{Let }$n$\textit{\ be a positive integer, }%
$n\geq 3$\textit{\ and let} $\xi $\textit{\ be a primitive root of order} $n$ of the unity.
 Let  $\mathbb{Q}\left(\xi_{n}\right)$ be the $n$th cyclotomic field. Let $p$ be a prime positive integer. Then  $p$ splits completely in the ring $\mathbb{Z}\left[\xi_{n}\right]$ 
if and only if $p\equiv 1$ (mod $n$).

\end{corollary}

\begin{corollary}
\label{2.5}
(\cite{ireland}) 
\textit{Let }$n$\textit{\ be a positive integer, }%
$n\geq 3$\textit{\ and let} $\xi $\textit{\ be a primitive root of order} $n$ of the unity.
 Let  $\mathbb{Q}\left(\xi_{n}\right)$ be the $n$th cyclotomic field. Let $p$ be a prime positive integer and let $P$ be a prime ideal in $\mathbb{Z}\left[\xi_{n}\right]$  such that $P\cap \mathbb{Z}=p \mathbb{Z}.$ If $p$ is odd then $P$ is ramified if and only if $p|n.$ If $p=2$
then $P$ is ramified if and only if $4|n.$

\end{corollary}

We now recall other important things about cyclotomic fields, which we need in our proofs.

\begin{proposition}
\label{propr.}
 (\cite{Washington}
 \textit{Let }$n$\textit{\ be a positive integer, }%
$n\geq 3$\textit{\ and let} $\xi_{n}$\textit{\ be a primitive root of order} $n$ of the unity.
 Let  $\mathbb{Q}\left(\xi_{n}\right)$ be the $n$th cyclotomic field. The maximal real subfield of the cyclotomic field 
$\mathbb{Q}\left(\xi_{n}\right)$ is  $\mathbb{Q}\left(\xi_{n}+ \xi^{-1}_{n}\right)$ and the degree of field extension 
$\mathbb{Q}\subset \mathbb{Q}\left(\xi_{n}+ \xi^{-1}_{n}\right)$ is $\frac{\varphi\left(n\right)}{2}$.
\end{proposition}

\begin{proposition}
\label{twodotseven}
 (\cite{Washington}).
Let $p$ be an odd prime  integer. Then $\mathbb{Q}\left(\sqrt{\left(-1\right)^{\frac{p-1}{2}}\cdot p}\right)$  is a quadratic subfield
of the  cyclotomic field $\mathbb{Q}\left(\xi_{p}\right).$
\end{proposition}
\begin{proposition}
\label{twodoteight}
(\cite{Savin}) Let $r$ and $s$ be two positive integers, with g.c.d.$\left(r,s\right)=1.$ Then $\mathbb{Q}\left(\xi_{r}, \xi_{s} \right)=\mathbb{Q}\left( \xi_{r\cdot s} \right).$
\end{proposition}

\begin{proposition}
\label{twodoteleven}
(\cite{Cox}) 
Let $K\subset  F \subset  L$ be two field extensions and let  $\mathcal{O}_{K},$  $\mathcal{O}_{F},$  $\mathcal{O}_{L}$ be the rings of algebraic integers.
A prime $p$ of  $\mathcal{O}_{K}$ splits completely in $\mathcal{O}_{L}$ if and only if $p$
splits completely in $\mathcal{O}_{F}$ and all primes of  $\mathcal{O}_{F}$ containing $p$ split completely in  $\mathcal{O}_{L}.$
\end{proposition}

\begin{proposition}
\label{*}
(\cite{janusz}) 
Let $K\subset  L$ be a Galois fields extension and let $\mathcal{O}_{K},$  $\mathcal{O}_{L}$ be the rings of algebraic integers. Let $p$ be a prime in  $\mathcal{O}_{K}.$
Since $\mathcal{O}_{L}$ is a Dedekind ring, the ideal $p\mathcal{O}_{L}$
is written uniquely in the form
$$p\mathcal{O}_{L}=P^{e_{1}}_{1}\cdot P^{e_{2}}_{2}\cdot...\cdot P^{e_{g}}_{g},$$
where $g, e_{1}, e_{2},...,e_{g} $ are positive integers and $P_{1}, P_{2},...., P_{g}$  are distinct prime ideals of the ring $\mathcal{O}_{L}$. Let the
completions fields $K_{p}$ and $L_{P}.$ Then:

 \begin{enumerate}[\rm1.]
 \item
$e_{1}=e_{2}=...=e_{g}=e$;
 \item
$\left[L : K\right]=efg,$ where  $f$ is the residual degree of each $P_{i},$ $i=\overline{1,g}$  in the fields extension $K\subset  L$, $f=\left[\mathcal{O}_{L}/P_{i} : \mathcal{O}_{K}/p\mathcal{O}_{k}\right]$;
\item
if $P$$\in$$\left\{P_{1}, P_{2},...., P_{g}\right\},$  then $\left[L_{P}: K_{p}\right]=ef.$

 \end{enumerate}

\end{proposition}

Now, we recall some results about central simple algebras.\\
Let $K$ be a field and  let $A$ be a central simple algebra over $K.$
It is known that the dimension $n$ of $A$\ over the field $K$ is a square. The positive
integer $\sqrt{n}$ is called\textit{\ the degree} of the algebra $A$.\newline
$A$ is called\textit{\ a division algebra} if the equations $\alpha x=\beta,\, y\alpha=\beta,\forall \alpha, \beta\in A$, $\alpha\neq 0,$ have unique solutions.
If $A$ is a finite-dimensional algebra, then $A$ is a division algebra if and only if $A$
is without zero divisors ($x\neq 0,y\neq 0\Rightarrow xy\neq 0$).\\
$A$ is called \textit{split} by $K$ if $A$ is isomorphic with a
matrix algebra over $K.$\newline
Let $K\subset F$ be a fields extension. $A$ is called \textit{split} by $F$ and $F$ is called
a \textit{splitting field} for $A$ if $A\otimes _{K}F$
is a matrix algebra over $F.$\\
In \cite{GS} there is the following criterion to decide whether an algebra of quaternions splits or not by a field $K.$
\begin{proposition}
\label{Proposition2.9}
 (\cite{GS}). \textit{Let} $K$ \textit{ be a field. Then the quaternion algebra} $\mathbb{H}%
_{K}\left(\alpha, \beta\right) $ \textit{splits if and only if} the
conic $C\left(\alpha, \beta\right)$ : $\alpha x^{2} +  \beta y^{2} = z^{2}$
has a rational point over $K$ (i.e. if there are $x_{0},y_{0},z_{0}$$\in$$K$ such that $\alpha x^{2}_{0}+ \beta y^{2}_{0}= z^{2}_{0}$).

\end{proposition}
We recall now Minkovski-Hasse Theorem  for a quadratic form in three variables.
\smallskip\\
\textbf{Minkovski-Hasse Theorem.} (\cite{BS}) 
\textit{The quadratic  form with  nonzero  rational coefficients} $\alpha x^{2}+\beta y^{2}- z^{2}$  \textit{represents zero in the field of rational numbers if and only if for all primes} $p$ (\textit{and also for} $p=\infty$), we have

$$\left(\frac{\alpha,\beta}{p}\right)=1,$$
\textit{where} $\left(\frac{\alpha,\beta}{p}\right)$ \textit{is the Hilbert symbol in the} $p$-\textit{adic field} $\mathbb{Q}_{p}.$

\begin{corollary}
\label{cases 1} (\cite{BS})\\
i) If $p$ is not equal with $2$ or $\infty$ and $p$ does not enter into the factorizations of $\alpha$ and $\beta$ into prime powers (which means that $\alpha$ and $\beta$ are $p$-adic units), then the form $\alpha x^{2}+\beta y^{2}- z^{2}$ represents zero in the $p$-adic fields $\mathbb{Q}_{p}$ and thus for all such $p$ the Hilbert symbol $\left(\frac{\alpha,\beta}{p}\right)=1.$\\
ii) $\left(\frac{\alpha,\beta}{\infty}\right)=1,$ if $\alpha>0$ or $\beta>0,$\\
$\left(\frac{\alpha,\beta}{\infty}\right)=-1,$ if $\alpha<0$ and $\beta<0,$\\
where $\left(\frac{\alpha,\beta}{\infty}\right)$ \textit{is the Hilbert symbol in the field} $\mathbb{R}.$
\end{corollary}
\begin{corollary}
\label{cases 2} (\cite{BS})\\
The product of the Hilbert symbols in the $p$-adic fields satisfy
$$\prod_{p}\left(\frac{\alpha,\beta}{p}\right)=1,$$
\textit{where} $p$-\textit{runs through all prime numbers and the symbol} $\infty.$
\end{corollary}

We recall the following result about the  quaternion algebras :
\begin{proposition}
\label{twodotfive}
(\cite{ledet}).
Let   $K$   be a field with   char $K \neq 2$   and let  $\alpha, \beta \in K \backslash \{  0 \}$.
 Then the quaternion algebra $H_{K} ( \alpha, \beta )$   is either a split or a division algebra.
\end{proposition}
We recall the following remark, which can be found in \cite{lam}.
\begin{remark}
 \label{twodotsix} Let $L/K$ be a finite extension of fields of odd degree, and let $\alpha, \beta \in K \backslash \{0\}$.
Then the quaternion algebra $H_{L}(\alpha, \beta)$ splits if and only if $H_{K}(\alpha, \beta)$ splits.
\end{remark}

Let $K$ be a number field and let a quaternion algebra $H_{K}\left(\alpha, \beta\right).$
 Let $\mathcal{O}_{K}$ be the ring of algebraic integers of the field $K.$ 
The reduced discriminant $D_{H_{K}\left(\alpha, \beta\right)}$ of the quaternion algebra $H_{K}\left(\alpha, \beta\right)$
is the integral ideal of the ring  $\mathcal{O}_{K}$ equal to the product of prime ideals of $\mathcal{O}_{K}$ that ramify in $H_{K}\left(\alpha, \beta\right)$ (see \cite{alsina}).\\
There is the following criterion to decide whether quaternion algebras split.
\begin{proposition}
\label{Proposition2.11}
  (\cite{alsina}). \textit{Let} $K$ \textit{ be a field. Then the quaternion algebra} $\mathbb{H}%
_{K}\left(\alpha, \beta\right) $ \textit{is split if and only if} $D_{H_{K}\left(\alpha, \beta\right)} = \mathcal{O}_{K}.$
\end{proposition}
In the case when  $\mathcal{O}_{K}$ is a principal ring, the ideals of ring  $\mathcal{O}_{K}$ identify with its generators, up to units.
So, a quaternion algebra $H_{\mathbb{Q}}\left(\alpha, \beta\right)$ is split if and only if $D_{H}=1.$ According to Proposition \ref{twodotfive},  a quaternion algebra $H_{\mathbb{Q}}\left(\alpha, \beta\right)$ is
a division algebra if and only if there is at least a prime positive integer $p$ which divides $D_{H\left(\alpha, \beta\right)}.$ According to a result from  \cite{lKohel}, if a prime $p$ $|$ $D_{H\left(\alpha, \beta\right)}$\ then $p$$|$$2\alpha\beta$.\ \ (2.1) \\ 
\smallskip\\
We recall the  following lemma about the discriminant of a quaternion algebra $H_{\mathbb{Q}}(p,q)$ .

\begin{lemma}
\label{twodottwelve}
 (\cite{alsina})
Let  $p$ and $q$ be two primes, and let $H_{\mathbb{Q}}(p,q)$ be a quaternion algebra of discriminant $D_H$.
\begin{enumerate}[\rm1.]
\item
If $p\equiv q\equiv 3$ $\pmod4$   {and} $(\frac{q}{p})\neq 1$,   {then} $D_{H}=2p$;
 \item
If $q=2$ and  $p\equiv 3$ $\pmod8$,  {then} $D_{H}=pq=2p$;
\item
If  $p$   {or} $q\equiv1 \pmod 4$, with  $p\neq q$  {and} $(\frac{p}{q})=-1$,   then  $D_{H}=pq$.
\end{enumerate}
\end{lemma}
We recall certain results of Albert-Brauer-Hasse-Noether, about the connection between a quaternion algebra  $H_{F}$ over a number field $F$ and the quaternion algebra  $H_{K}$ over a number field $K,$ when $K$ is a quadratic extension of $F.$ 

\begin{theorem}
\label{Theorem2.13}(\cite{llinowitz}) (Albert-Brauer-Hasse-Noether). \textit{Let} $H_{F}$ \textit{be a quaternion algebra over a number field} $F$ \textit{and let} $K$ \textit{be a quadratic field extension of} $F.$ \textit{Then there is an embedding of} $K$ \textit{into}
$H_{F}$ \textit{if and only if no prime of} $F$ \textit{which ramifies in} $H_{F}$ \textit{splits in} $K.$
\end{theorem}

\begin{proposition}
\label{Proposition2.14}
(\cite{Voight}). \textit{Let} $F$ \textit{be a number field and let} $K$ \textit{be a field containing} $F.$ \textit{Let} $H_{F}$ \textit{be a quaternion algebra over} $F.$ \textit{Let} $H_{K} =H_{F}\otimes _{F} K $ \textit{be a quaternion algebra over
} $K.$ \textit{If} $[K : F] = 2,$ \textit{then} $K$ \textit{splits} $H_{F}$ \textit{if and
only if there exists an} $F$-\textit{embedding} $K\hookrightarrow H_{F}.$
\end{proposition}

\begin{proposition}
\label{Proposition2.15}
(\cite{Voight1}). \textit{Let} $H_{F}\left(a,b\right)$ be a division quaternion algebra over a local field $F,$ and let $L$ be a separable field extension of $F$ of finite degree. Then $L$ is a spliting field for $H_{F}\left(a,b\right)$ (i.e $H_{L}\left(a,b\right)$ splits)
if and only if the degree $\left[L:F\right]$ is even.

\end{proposition}

\begin{equation*}
\end{equation*}
\section{ Division quaternion algebras over a cyclotomic field  $\mathbb{Q}\left(\xi_{n}\right),$  with $n$$\in$$\left\{3,4,5,6,7,8,9,10,11,12\right\}$} 
\label{three}
\begin{equation*}
\end{equation*}
In the paper \cite{astz} (Theorem 4.6), V. Acciaro, D. Savin, M. Taous and A. Zekhnini found necessary and sufficient conditions for that quaternion algebras  $H_{\mathbb{Q}(\sqrt d)}(p, q)$ to be division algebras over quadratic fields  $\mathbb{Q}(\sqrt d).$

\begin{theorem}
\label{theorem3.1}
(Theorem 4.6 from  \cite{astz}). {Let} $d\neq1$  {be a squarefree integer} and let
 $K=\mathbb{Q}(\sqrt{d})$,
 with discriminant  $\Delta_{K}$. {Let} $p$ and $q$  be two  positive primes.
 {Then the quaternion algebra} $H_{K}(p,q)$   {is a division algebra if and only if} one of the following conditions
  holds:
 \begin{enumerate}[\rm1.]
 \item
$p$ and $q$  are odd and distinct, $p$ or $q \equiv$ $1$ (mod $4$),
$(\frac{p}{q})=-1$
     and
    $(\frac{\Delta_{K}}{p})=1$   {or}  $(\frac{\Delta_{K}}{q})=1$;
 \item
  $q=2$,  $p \equiv$ $3$  (mod $8$),   or      $p \equiv$ $5$ (mod $8$) and
 either $(\frac{\Delta_{K}}{p})=1$   {or}
  $d\equiv$$1$ (mod $8$);
 \item
 $p$ and $q$ are odd,  with $p\equiv q\equiv 3$ (mod $4$), and
 \begin{itemize}
 \item
  $(\frac{q}{p})\neq 1$ and
 either $(\frac{\Delta_{K}}{p})=1$   {or} $d$$\equiv$$1$ (mod $8$);\\
 or
 \item
  $(\frac{p}{q})\neq 1$ and
 either $(\frac{\Delta_{K}}{q})=1$   {or} $d$$\equiv$$1$ (mod $8$).
 \end{itemize}
 \end{enumerate}
\end{theorem}

We make the following remark.
\begin{remark}
\label{remark3.2}
 Let $L/K$ be a finite extension of fields of odd degree, and let $\alpha, \beta \in K \backslash \{0\}$.
Then the quaternion algebra $H_{L}(\alpha, \beta)$ is a division algebra if and only if $H_{K}(\alpha, \beta)$ is a division algebra.
\end{remark}

\begin{proof}
The proof follows immediately, using Remark \ref{twodotsix} and Proposition\ref{twodotfive}.
\end{proof}

We study now quaternion algebras $H_{\mathbb{Q}\left(\xi_{n}\right)}\left(p_{1}, p_{2}\right)$  over a cyclotomic field  $\mathbb{Q}\left(\xi_{n}\right),$  where $n$$\in$$\left\{3,4,6,7,8,9,11,12\right\}$ More precisely, we find necessary and sufficient conditions for such an algebra to be a division algebra.
\\

For $n=3,$ the 3rd cyclotomic field  is $\mathbb{Q}\left(\xi_{3}\right)=$ $\mathbb{Q}(i\sqrt 3).$ Necessary and sufficient conditions for quaternion algebra   $H_{\mathbb{Q}\left(\xi_{3}\right)}\left(p_{1}, p_{2}\right)$  to be a division algebra result from the articles  \cite{astz} (Theorem 4.6) and  \cite{astz1}.
We will give all the details later (in the Proposition \ref{threedotsix})\\
For $n=6,$ the 6th cyclotomic field, applying Proposition \ref{twodoteight}, we obtain $\mathbb{Q}\left(\xi_{6}\right)=\mathbb{Q}(\xi_{2}, \xi_{3})=\mathbb{Q}(-1, i\sqrt 3)=\mathbb{Q}(i\sqrt 3),$ so the case $n=6$ is reduced to the case $n=3.$\\
For $n=4,$ the 4th cyclotomic field  is  $\mathbb{Q}\left(\xi_{4}\right)=$ $\mathbb{Q}(i).$ Necessary and sufficient conditions for quaternion algebra    $H_{\mathbb{Q}\left(i\right)}\left(p_{1}, p_{2}\right)$ to be a division algebra result from the article  \cite{astz} (Theorem 4.6).\\
For $n=7,$ we obtain the following result:
\begin{proposition}
\label{threedotthree} Let $\xi_{7}$ be a primitive root of order $7$ of the unity and  let  $\mathbb{Q}\left(\xi_{7}\right)$ be the $7$th cyclotomic field. Let $p_{1}, p_{2}$ be two distinct prime integers.
Then the quaternion algebra $H_{\mathbb{Q}\left(\xi_{7}\right)}\left(p_{1}, p_{2}\right)$ is a division algebra if and only if one of the following conditions
  holds:
 \begin{enumerate}[\rm1.]
 \item
$p_{1}$ and $p_{2}$  are odd and distinct, $p_{1}$ or $p_{2} \equiv 1$ (mod $4$),
$(\frac{p_{1}}{p_{2}})=-1$
     and
    $(\frac{-7}{p_{1}})=1$   {or}  $(\frac{-7}{p_{2}})=1$;
 \item
  $p_{2}=2$,  $p_{1} \equiv 3$ (mod $8$)   or      $p_{1} \equiv 5$ (mod $8$) ;

 \item
 $p_{1}$ and $p_{2}$ are odd,  with $p_{1}\equiv p_{2}\equiv 3$ (mod $4$).

 \end{enumerate}

\end{proposition}
\begin{proof}
According to Proposition   \ref{twodotseven}  we have the following field extensions
 $$\mathbb{Q}\subset \mathbb{Q}\left(i\sqrt{7}\right)\subset \mathbb{Q}\left(\xi_{7}\right).$$ 
We know the degrees $ \left[\mathbb{Q}\left(\xi_{7}\right): \mathbb{Q}\right]=\varphi\left(7\right)=6$
and $\left[\mathbb{Q}\left(i\sqrt{7}\right): \mathbb{Q}\right]=2.$ According to the transitivity theorem of finite extensions, it follows that $ \left[\mathbb{Q}\left(\xi_{7}\right): \mathbb{Q}\left(i\sqrt{7}\right)\right]=3$.\\
Applying  Remark \ref{remark3.2}, Theorem \ref{theorem3.1} (for $d=-7\equiv 1$ (mod $8$)) and  Proposition \ref{twodotfthree} (4.), we obtain that  the quaternion algebra $H_{\mathbb{Q}\left(\xi_{7}\right)}\left(p_{1}, p_{2}\right)$ is a division algebra if and only if one of the conditions 1., 2., 3 from the conclusion of Proposition \ref{threedotthree}.
is fulfilled.
\end{proof}

In the paper \cite{astz2} (Theorem 4.6), V. Acciaro, D. Savin, M. Taous and A. Zekhnini found necessary and sufficient conditions for that quaternion algebras  $H_{\mathbb{Q}(\sqrt d_{1}, \sqrt d_{2}}(p, q)$ to be division algebras over quadratic fields  $\mathbb{Q}(\sqrt d_{1}, \sqrt d_{2}).$

\begin{theorem}
\label{theorem3.4} 
(Theorem 4.6 from  \cite{astz2})
{Let} $d_{1}$  {and} $d_{2}$  {be distinct   squarefree integers} not equal to one. {Let} $K=\mathbb{Q}(\sqrt{d_{1}}, \sqrt{d_{2}})$, {and let} $K_{i}=\mathbb{Q}(\sqrt{d_{i}})$ ($i=1,2$) with discriminant  $\Delta_{K_{i}}$. {Let} $p$ and $q$  be two positive
primes. 
 {Then the quaternion algebra} $H_{K}(p,q)$   {is a division algebra if and only if} one of the following conditions
  holds:
 \begin{enumerate}
 \item
$p$ and $q$  are odd and distinct, and \\
$(\frac{p}{q})=-1$,   and \\
     $p \equiv 1 $ ({mod} $4$)  {or} $q  \equiv 1 $ ({mod} $4$),  and  \\
    $(\frac{\Delta_{K_{1}}}{p})=(\frac{\Delta_{K_{2}}}{p})=1$   {or}  $(\frac{\Delta_{K_{1}}}{q})=(\frac{\Delta_{K_{2}}}{q})=1$;
 \item
  $q=2$, and  \\
 $p \equiv 3$ ({mod} $8$)   or      $p \equiv 5$ ({mod} $8$), and \\
 either $(\frac{\Delta_{K_{1}}}{p})=(\frac{\Delta_{K_{2}}}{p})=1$   {or} 
  $d_{1}, d_{2}$$\equiv$$1$ ({mod} $8$);
 \item
 $p$ and $q$ are odd,  with $p\equiv q\equiv 3$ ({mod} $4$), and
 \begin{itemize} 
 \item
  $(\frac{q}{p})\neq 1$, and
  \item
 either $(\frac{\Delta_{K_{1}}}{p})=(\frac{\Delta_{K_{2}}}{p})=1$   {or} $d_{1}, d_{2}$$\equiv$$1$ ({mod} $8$);
 \end{itemize}
 or
 \begin{itemize} 
 \item
  $(\frac{p}{q})\neq 1$, and
  \item
 either $(\frac{\Delta_{K_{1}}}{q})=(\frac{\Delta_{K_{2}}}{q})=1$   {or} $d_{1}, d_{2}$$\equiv$$1$ ({mod} $8$);
 \end{itemize}
 \end{enumerate} 
\end{theorem}

For the case when $n=8,$ we obtain the following result:
\begin{proposition}
\label{threedotfive} Let $\xi_{8}$ be a primitive root of order $8$ of the unity and  let  $\mathbb{Q}\left(\xi_{8}\right)$ be the $8$th cyclotomic field. Let $p_{1}, p_{2}$ be two distinct prime integers.
Then the quaternion algebra $H_{\mathbb{Q}\left(\xi_{8}\right)}\left(p_{1}, p_{2}\right)$ is a division algebra if and only if 

$$p_{1}\: and \: p_{2} \:are\:  odd\:  and\:  distinct,\:  and\:  (\frac{p_{1}}{p_{2}})=-1, \:  and $$
     $$p_{1} \equiv 1\:  ({mod}\: 8) \: {or}\: p_{2}\:  \equiv 1\:  ({mod}\: 8). $$
 
\end{proposition}
\begin{proof}
The order of the Galois group Gal$\left(\mathbb{Q}\left(\xi_{8}\right)/\mathbb{Q}\right)$ is $\varphi\left(8\right)=4$ and it is isomorphic to the group $\left(\mathbb{Z}_{2}\times \mathbb{Z}_{2}, +\right).$\\
We have the following field extensions $\mathbb{Q}\subset \mathbb{Q}\left(\xi_{8}+\xi^{-1}_{8}\right)= \mathbb{Q}\left(cos\frac{2\pi}{8}\right)  \subset \mathbb{Q}\left(\xi_{8}\right).$\\
$cos\frac{2\pi}{8}=cos\frac{\pi}{4}=\frac{\sqrt{2}}{2},$ so $\mathbb{Q}\left(\xi_{8}+\xi^{-1}_{8}\right)=\mathbb{Q}\left(\sqrt{2}\right).$\\
$\mathbb{Q}\left(\xi_{8}\right)= \mathbb{Q}\left(cos\frac{\pi}{4}+i\cdot sin\frac{\pi}{4}\right)=\mathbb{Q}\left(\frac{\sqrt{2}}{2}+i\cdot \frac{\sqrt{2}}{2}\right)=\mathbb{Q}\left(i, \sqrt{2}\right).$\\
So, the previous field extensions become $\mathbb{Q}\subset \mathbb{Q}\left(\sqrt{2}\right) \subset \mathbb{Q}\left(i, \sqrt{2}\right)=\mathbb{Q}\left(\xi_{8}\right).$  We have also the following field extensions
$\mathbb{Q}\subset \mathbb{Q}\left(i\right) \subset \mathbb{Q}\left(i, \sqrt{2}\right)=\mathbb{Q}\left(\xi_{8}\right).$ From this point, we can apply Theorem \ref{theorem3.4}.\\
Noting $K_1=\mathbb{Q}\left(i\right)$ and $K_2= \mathbb{Q}\left(\sqrt{2}\right),$
we obtain the Legendre's symbols $(\frac{\Delta_{K_{1}}}{p_{1}})=(\frac{-1}{p_{1}}),$  $(\frac{\Delta_{K_{2}}}{p_{1}})=(\frac{2}{p_{1}}).$\\
If  $p_{2}=2,$ $p_{1} \equiv 3$ (mod $8$)   or $p_{1} \equiv 5$ (mod $8$), it results 
$(\frac{2}{p_{1}})=-1.$ Also $-1\not\equiv 1$ (mod $8$). It follows that we cannot have the case (ii) from Theorem \ref{theorem3.4}. Similarly, we obtain that we cannot have the case (iii) from Theorem \ref{theorem3.4}.\\
Applying Theorem  \ref{theorem3.4}, it results that $p_{1}$ and $p_{2}$ are odd  and distinct, and $(\frac{p_{1}}{p_{2}})=-1$, and 
$p_{1} \equiv 1$ ({mod} $4$)  or $p_{2}\equiv 1$ (mod $4$),  and  
    $(\frac{-1}{p_{1}})=(\frac{2}{p_{1}})=1$  or $ (\frac{-1}{p_{2}})=(\frac{2}{p_{2}})=1.$  Now applying Proposition \ref{twodotfthree}, we immediately obtain that the conditions 
$p_{1} \equiv 1$ ({mod} $4$)  or $p_{2}\equiv 1$ (mod $4$),  and  
    $(\frac{-1}{p_{1}})=(\frac{2}{p_{1}})=1$  or $ (\frac{-1}{p_{2}})=(\frac{2}{p_{2}})=1$ are equivalent to $p_{1} \equiv 1$ ({mod} $8$)  or $p_{2}\equiv 1$ (mod $8$).
\end{proof}
For the case when $n=9,$ we obtain the following result:
\begin{proposition}
\label{threedotsix}
 Let $\xi_{9}$ be a primitive root of order $9$ of the unity and  let  $\mathbb{Q}\left(\xi_{9}\right)$ be the $9$th cyclotomic field. Let $p_{1}, p_{2}$ be two distinct prime integers.
Then the quaternion algebra $H_{\mathbb{Q}\left(\xi_{9}\right)}\left(p_{1}, p_{2}\right)$ is a division algebra if and only if one of the following conditions
  holds:
 \begin{enumerate}[\rm1.]
 \item
 $p_{1}$ and $p_{2}$  are odd and distinct, $p_{1}$ or $p_{2} \equiv$ $1$ (mod $12$), $(\frac{p_{1}}{p_{2}})=-1;$
 \item
$p_{2}=2$,  $p_{1} \equiv$ $19$  (mod $24$),   or   $p_{1} \equiv$ $13$ (mod $24$);

 \item
 $p_{1}$ and $p_{2}$ are odd,  with $p_{1}\equiv p_{2}\equiv 3$ (mod $4$), and
 \begin{itemize}
 \item
  $(\frac{p_{2}}{p_{1}})\neq 1$ and $p_{1} \equiv$ $1$  (mod $3$);\\
 or
 \item
  $(\frac{p_{1}}{p_{2}})\neq 1$ and $p_{2} \equiv$ $1$  (mod $3$).
 \end{itemize}
 \end{enumerate}

\end{proposition}
\begin{proof}
$\left|Gal\left(\mathbb{Q}\left(\xi_{9}\right)/\mathbb{Q}\right)\right|=\varphi\left(9\right)=6$ and the Galois group $Gal\left(\mathbb{Q}\left(\xi_{9}\right)/\mathbb{Q}\right)$ is isomorphic to the group $\left(\mathbb{Z}_{6}, +\right),$ so 
$Gal\left(\mathbb{Q}\left(\xi_{9}\right)/\mathbb{Q}\right)$  is a cyclic group, $Gal\left(\mathbb{Q}\left(\xi_{9}\right)/\mathbb{Q}\right)=<\tau>,$ where $\tau\left(\xi_{9}\right)=\xi^{2}_{9}.$ 
The Galois group  Gal$\left(\mathbb{Q}\left(\xi_{9}\right)/\mathbb{Q}\right)$ contains the following two proper subgroups: $H_{1}=\left\{1_{\mathbb{Q}\left(\xi_{9}\right)}, \tau^{3}\right)$ and 
$H_{2}=\left\{1_{\mathbb{Q}\left(\xi_{9}\right)}, \tau^{2}, \tau^{4}\right).$ 
From the Galois correspondence, it results immediately that the proper subfields of the 9th cyclotomic field $\mathbb{Q}\left(\xi_{9}\right)$ are:
$\mathbb{Q}\left(\xi_{9}\right)^{H_{1}}=\mathbb{Q}\left(\xi_{9}+\xi^{-1}_{9}\right)= \mathbb{Q}\left(cos\frac{2\pi}{9}\right)$ and $\mathbb{Q}\left(\xi_{9}\right)^{H_{2}}=\mathbb{Q}\left(\xi^{3}_{9}\right)= \mathbb{Q}\left(\xi_{3}\right).$
We have the following field extensions 
$$\mathbb{Q}\subset \mathbb{Q}\left(\xi_{3}\right)=\mathbb{Q}\left(i\sqrt{3}\right) \subset \mathbb{Q}\left(\xi_{9}\right).$$
The degrees of the following field extensions are known: 
$$\left[\mathbb{Q}\left(\xi_{9}\right): \mathbb{Q}\right]=6, \left[\mathbb{Q}\left(i\sqrt{3}\right): \mathbb{Q}\right]=2.$$ 
It results that $\left[\mathbb{Q}\left(\xi_{9}\right): \mathbb{Q}\left(i\sqrt{3}\right)\right]=3,$ 
According to Remark \ref{remark3.2}, $H_{\mathbb{Q}\left(\xi_{9}\right)}\left(p_{1}, p_{2}\right)$ is a division algebra if and only if $H_{\mathbb{Q}\left(i\sqrt{3}\right)}\left(p_{1}, p_{2}\right)$ is a division algebra.
Applying  Theorem \ref{theorem3.1} when $d=-3$ , we obtain that  the quaternion algebra $H_{\mathbb{Q}\left(\xi_{9}\right)}\left(p_{1}, p_{2}\right)$ is a division algebra if and only if
\begin{enumerate}[\rm1.]
 \item
$p_{1}$ and $p_{2}$  are odd and distinct, $p_{1}$ or $p_{2} \equiv$ $1$ (mod $4$),
$(\frac{p_{1}}{p_{2}})=-1$
     and
    $(\frac{-3}{p_{1}})=1$   {or}  $(\frac{-3}{p_{2}})=1$;
 \item
  $p_{2}=2$,  $p_{1} \equiv$ $3$  (mod $8$),   or     $p_{1} \equiv$ $5$ (mod $8$) and
 $(\frac{-3}{p_{1}})=1$;  
 \item
 $p_{1}$ and $p_{2}$ are odd,  with $p_{1}\equiv p_{2}\equiv 3$ (mod $4$), and
 \begin{itemize}
 \item
  $(\frac{p_{2}}{p_{1}})\neq 1$ and $(\frac{-3}{p_{1}})=1$;\\
 or
 \item
  $(\frac{p_{1}}{p_{2}})\neq 1$ and $(\frac{-3}{p_{2}})=1$.
 \end{itemize}
 \end{enumerate}
Applying Proposition \ref{twodotfthree} and the Chinese remainder Lemma, condition 1 is equivalent to $p_{1}$ and $p_{2}$  are odd and distinct, $p_{1}$ or $p_{2} \equiv$ $1$ (mod $12$), $(\frac{p_{1}}{p_{2}})=-1$;\\
condition 2. is equivalent to  $p_{2}=2$,  $p_{1} \equiv$ $19$  (mod $24$),   or     $p_{1} \equiv$ $13$ (mod $24$);\\
condition 3. is equivalent to 
 $p_{1}$ and $p_{2}$ are odd,  with $p_{1}\equiv p_{2}\equiv 3$ (mod $4$), and
 \begin{itemize}
 \item
  $(\frac{p_{2}}{p_{1}})\neq 1$ and $p_{1} \equiv$ $1$  (mod $3$);\\
 or
 \item
  $(\frac{p_{1}}{p_{2}})\neq 1$ and  $p_{2} \equiv$ $1$  (mod $3$).
 \end{itemize}
So,  the quaternion algebra $H_{\mathbb{Q}\left(\xi_{9}\right)}\left(p_{1}, p_{2}\right)$ is a division algebra if and only if one of the conditions 1., 2., 3 from the conclusion of Proposition \ref{threedotsix}
is fulfilled.
\end{proof}
For the case when $n=11,$ we obtain the following result:
\begin{proposition}
\label{threedotseven}
Let $\xi_{11}$ be a primitive root of order $11$ of the unity and  let  $\mathbb{Q}\left(\xi_{11}\right)$ be the $11$th cyclotomic field. Let $p_{1}, p_{2}$ be two distinct prime integers.
Then the quaternion algebra $H_{\mathbb{Q}\left(\xi_{11}\right)}\left(p_{1}, p_{2}\right)$ is a division algebra if and only if one of the following conditions
  holds:
\begin{enumerate}[\rm1.]
 \item
$p_{1}$ and $p_{2}$  are odd and distinct, $p_{1}$ or $p_{2} \equiv 1$ (mod $4$),
$(\frac{p_{1}}{p_{2}})=-1$
     and
    $(\frac{-11}{p_{1}})=1$   {or}  $(\frac{-11}{p_{2}})=1$;
 \item
  $p_{2}=2$,  $p_{1} \equiv 3$ (mod $8$)   or      $p_{1} \equiv 5$ (mod $8$) and $(\frac{-11}{p_{1}})=1$;

 \item
 $p_{1}$ and $p_{2}$ are odd,  with $p_{1}\equiv p_{2}\equiv 3$ (mod $4$) and
 \begin{itemize}
 \item
  $(\frac{p_{2}}{p_{1}})\neq 1$ and $(\frac{-11}{p_{1}})=1$;\\
 or
 \item
  $(\frac{p_{1}}{p_{2}})\neq 1$ and $(\frac{-11}{p_{2}})=1$.
 \end{itemize}

 \end{enumerate}

\end{proposition}

\begin{proof}
According to Proposition \ref{twodotseven}  we have the following field extensions
$$\mathbb{Q}\subset \mathbb{Q}\left(i\sqrt{11}\right) \subset \mathbb{Q}\left(\xi_{11}\right).$$
Since the degrees $ \left[\mathbb{Q}\left(\xi_{11}\right): \mathbb{Q}\right]=\varphi\left(11\right)=10$
and $\left[\mathbb{Q}\left(i\sqrt{11}\right): \mathbb{Q}\right]=2,$ 
 it follows that the degree $\left[\mathbb{Q}\left(\xi_{11}\right)  : \mathbb{Q}\left(i\sqrt{11}\right)\right]=5.$ Applying  Remark \ref{remark3.2} and Theorem \ref{theorem3.1} (for $d= -11$), we obtain that  the quaternion algebra $H_{\mathbb{Q}\left(\xi_{11}\right)}\left(p_{1}, p_{2}\right)$ is a division algebra if and only if one of the conditions 1., 2., 3. from the conclusion of Proposition \ref{threedotseven}.
is fulfilled.
\end{proof}
For the case when $n=12,$ we obtain the following result:
\begin{proposition}
\label{threedoteight}
Let $\xi_{12}$ be a primitive root of order $12$ of the unity and  let  $\mathbb{Q}\left(\xi_{12}\right)$ be the $12$th cyclotomic field. Let $p_{1}, p_{2}$ be two distinct prime integers.
Then the quaternion algebra $H_{\mathbb{Q}\left(\xi_{12}\right)}\left(p_{1}, p_{2}\right)$ is a division algebra if and only if one of the following conditions
  holds:
\begin{enumerate}
 \item
$p_{1}$ and $p_{2}$  are odd and distinct, and \\
$(\frac{p_{1}}{p_{2}})=-1$,   and \\
     $p_{1} \equiv 1 $ ({mod} $12$)  {or} $p_{2}  \equiv 1 $ ({mod} $12$);
 \item
  $p_{2}=2$, and  $p_{1} \equiv 13$ ({mod} $24$).

 \end{enumerate}

\end{proposition}
\begin{proof}
According to Proposition \ref{twodoteight}, $\mathbb{Q}\left( \xi_{12} \right)=\mathbb{Q}\left(\xi_{4}, \xi_{3} \right)=\mathbb{Q}\left(i, i\sqrt{3} \right).$\\
Noting $K_1=\mathbb{Q}\left(i\right)$ and $K_2= \mathbb{Q}\left(i\sqrt{3}\right),$
we obtain the Legendre's symbols $(\frac{\Delta_{K_{1}}}{p_{1}})=(\frac{-1}{p_{1}}),$  $(\frac{\Delta_{K_{2}}}{p_{1}})=(\frac{-3}{p_{1}}).$\\
Applying  Theorem \ref{theorem3.4}, we obtain that  the quaternion algebra $H_{\mathbb{Q}\left(\xi_{12}\right)}\left(p_{1}, p_{2}\right)$ is a division algebra if and only if one of the following conditions
  holds:
\begin{enumerate}
 \item
$p_{1}$ and $p_{2}$  are odd and distinct, and \\
$(\frac{p_{1}}{p_{2}})=-1$,   and \\
     $p_{1} \equiv 1 $ ({mod} $4$)  {or} $p_{2}  \equiv 1 $ ({mod} $4$),  and  \\
    $(\frac{-1}{p_{1}})=(\frac{-3}{p_{1}})=1$   {or}  $(\frac{-1}{p_{2}})=(\frac{\-3}{p_{2}})=1$;
 \item
  $p_{2}=2$, and  \\
 $p_{1} \equiv 3$ ({mod} $8$)   or      $p_{1} \equiv 5$ ({mod} $8$), and \\
 either $(\frac{-1}{p_{1}})=(\frac{-3}{p_{1}})=1$   {or} 
  $d_{1}, d_{2}$$\equiv$$1$ ({mod} $8$);
 \item
 $p_{1}$ and $p_{2}$ are odd,  with $p_{1}\equiv p_{2}\equiv 3$ ({mod} $4$), and
 \begin{itemize} 
 \item
  $(\frac{p_{2}}{p_{1}})\neq 1$, and
  \item
 either $(\frac{-1}{p_{1}})=(\frac{-3}{p_{1}})=1$   {or} $d_{1}, d_{2}$$\equiv$$1$ ({mod} $8$);
 \end{itemize}
 or
 \begin{itemize} 
 \item
  $(\frac{p_{1}}{p_{2}})\neq 1$, and
  \item
 either  $(\frac{-1}{p_{2}})=(\frac{-3}{p_{2}})=1$   {or} $d_{1}, d_{2}$$\equiv$$1$ (mod $8$);
 \end{itemize}

 \end{enumerate}

Applying Proposition \ref{twodotfthree} and the Chinese remainder Lemma, condition (i) is 
equivalent to $p_{1}$ and $p_{2}$  are odd and distinct, 
$(\frac{p_{1}}{p_{2}})=-1$, $p_{1} \equiv 1 $ ({mod} $12$)  {or} $p_{2}  \equiv 1 $ ({mod} $12$);
condition (ii) is equivalent to  $p_{2}=2$, and   $p_{1} \equiv 5$ ({mod} $8$), and $(\frac{-1}{p_{1}})=(\frac{-3}{p_{1}})=1,$ which is
equivalent to  $p_{2}=2$ and   $p_{1} \equiv 13$ ({mod} $24$).\\
If $p_{1}\equiv p_{2}\equiv 3$ ({mod} $4$), it results $(\frac{-1}{p_{1}})=-1$ and $(\frac{-1}{p_{2}})=-1.$ Also $d_{1}=-1$$\not\equiv$$1$ (mod $8$),  $d_{2}=-3$$\not\equiv$$1$ (mod $8$), 
we cannot have condition(iii) from Theorem \ref{theorem3.4}.\\
So,  the quaternion algebra $H_{\mathbb{Q}\left(\xi_{12}\right)}\left(p_{1}, p_{2}\right)$ is a division algebra if and only if one of the conditions (i), (ii)  from the conclusion of Proposition \ref{threedoteight}
is fulfilled.

\end{proof}

We are studying quaternion algebras $H_{\mathbb{Q}\left(\xi_{n}\right)}\left(p_{1}, p_{2}\right)$  over a cyclotomic field  $\mathbb{Q}\left(\xi_{n}\right),$  where $n$$\in$$\left\{5,10\right\}$ More precisely, we give sufficient conditions for such an algebra to be a division algebra.\\
\begin{proposition}
\label{threedotnine}
Let $\xi_{5}$ be a primitive root of order $5$ of the unity and  let  $\mathbb{Q}\left(\xi_{5}\right)$ be the $5$th cyclotomic field. Let $p_{1}, p_{2}$ be two distinct odd prime integers, $p_{1}\equiv$$1$ (mod $5$), $\left(\frac{p_{2}}{p_{1}}\right)=-1.$ 
Then the quaternion algebra $H_{\mathbb{Q}\left(\xi_{5}\right)}\left(p_{1}, p_{2}\right)$ is a division algebra.

\end{proposition}

\begin{proof} According to Proposition \ref{twodotseven}, we have the following field extensions:
$$\mathbb{Q}\subset  \mathbb{Q}\left(\sqrt{5}\right)    \subset \mathbb{Q}\left(\xi_{5}\right),$$
with the degrees $\left[\mathbb{Q}\left(\sqrt{5}\right):\mathbb{Q}\right]=2$ and $\left[\mathbb{Q}\left(\xi_{5}\right):\mathbb{Q}\right]=\varphi\left(5\right) =4.$ It results that 
$\left[\mathbb{Q}\left(\xi_{5}\right): \mathbb{Q}\left(\sqrt{5}\right)\right]=2.$\\
If $p_{1}\equiv$$1$ (mod $5$),  applying Proposition \ref{twodotfthree} we have: $\left(\frac{p_{1}}{5}\right)=1$ and $\left(\frac{p_{1}}{5}\right)\cdot \left(\frac{5}{p_{1}}\right)=1$
it follow that $\left(\frac{5}{p_{1}}\right)=1$  According to Theorem \ref{theorem3.1}, the quaternion algebra $H_{\mathbb{Q}\left(\sqrt{5}\right)}\left(p_{1}, p_{2}\right)$ is a division algebra.
According to Lemma \ref{twodottwelve}, $p_{1} |$$D_{H\left(p_{1}, p_{2}\right)}$,  so $p_{1}$ is a prime which ramifies in this quaternion algebra.\\
The ring of algebraic integers of the quadratic field $\mathbb{Q}\left(\sqrt{5}\right)$ is $\mathbb{Z}\left[\frac{1+\sqrt{5}}{2}\right].$ \\
Since  $\left(\frac{5}{p_{1}}\right)=1$, according to Theorem \ref{theorem2.1}, $p_{1}$ splits completely in the ring  $\mathbb{Z}\left[\frac{1+\sqrt{5}}{2}\right].$ \\
Since $p_{1}\equiv$$1$ (mod $5$), according to Corollary \ref{twodotnine}, $p_{1}$ splits completely in the ring  $\mathbb{Z}\left[\xi_{5}\right].$ \\
Applying Proposition \ref{twodoteleven}, it follows that all primes of $\mathbb{Z}\left[\frac{1+\sqrt{5}}{2}\right]$ containing $p_{1}$ split completely in   $\mathbb{Z}\left[\xi_{5}\right].$
Applying Theorem \ref{Theorem2.13}, Proposition \ref{Proposition2.14} and Proposition \ref{twodotfive}, we obtain that the quaternion algebra $H_{\mathbb{Q}\left(\xi_{5}\right)}\left(p_{1}, p_{2}\right)$ is a division algebra.

\end{proof}

\begin{proposition}
\label{threedotten}
Let $\xi_{10}$ be a primitive root of order $10$ of the unity and  let  $\mathbb{Q}\left(\xi_{10}\right)$ be the $10$th cyclotomic field. Let $p_{1}, p_{2}$ be two distinct odd prime integers, $\left(\frac{p_{2}}{p_{1}}\right)=-1$ and $p_{1}\equiv$$1$ (mod $5$).
Then the quaternion algebra $H_{\mathbb{Q}\left(\xi_{10}\right)}\left(p_{1}, p_{2}\right)$ is a division algebra.

\end{proposition}

\begin{proof} 
According to Proposition \ref{twodoteight}, we have $\mathbb{Q}\left(\xi_{10}\right)=\mathbb{Q}\left(\xi_{5}, \xi_{2} \right)=\mathbb{Q}\left(\xi_{5}, -1 \right)= \mathbb{Q}\left(\xi_{5} \right).$
Applying Proposition \ref{threedotnine}, we obtain that the quaternion algebra $H_{\mathbb{Q}\left(\xi_{10}\right)}\left(p_{1}, p_{2}\right)$ is a division algebra.

\end{proof}

\begin{equation*}
\end{equation*}
\section{ Division quaternion algebras over a cyclotomic field  $\mathbb{Q}\left(\xi_{n}\right),$  when $n=l^{k},$  with $l$  is a prime integer
$l\equiv 3$ (mod $4$) } 
\label{four}
\begin{equation*}
\end{equation*}
Now we generalize Propositions  \ref{threedotthree}, \ref{threedotsix}, \ref{threedotseven} and  we obtain necessary and sufficient conditions for that the same quaternion algebras to be division algebras over the  $l^{k}$th cyclotomic field $\mathbb{Q}\left(\xi_{l^{k}}\right)$, where $l$ is an odd prime integer, $l \equiv 3$ ({mod} $4$) .

\begin{proposition}
\label{fourdotone}  Let $l$   {be an odd prime integer,} $l \equiv 3$ ({mod} $4$) and  let  $\mathbb{Q}\left(\xi_{l^{k}}\right)$ be the $l^{k}$th cyclotomic field, where $k$ is a positive integer. Let $p_{1}, p_{2}$ be two distinct prime integers.
Then the quaternion algebra $H_{\mathbb{Q}\left(\xi_{l^{k}}\right)}\left(p_{1}, p_{2}\right)$ is a division algebra if and only if one of the following conditions
  holds:
 \begin{enumerate}[\rm1.]
 \item
$p_{1}$ and $p_{2}$  are odd and distinct, $p_{1}$ or $p_{2} \equiv 1$ (mod $4$),
$(\frac{p_{1}}{p_{2}})=-1$
     and
    $(\frac{-l}{p_{1}})=1$   {or}  $(\frac{-l}{p_{2}})=1$;
 \item
  $p_{2}=2$,  $p_{1} \equiv 3$ (mod $8$)   or      $p_{1} \equiv 5$ (mod $8$) and
 either $(\frac{-l}{p_{1}})=1$   {or}
  $l\equiv7\pmod8$;
 \item
 $p_{1}$ and $p_{2}$ are odd,  with $p_{1}\equiv p_{2}\equiv 3$ (mod $4$), and
 \begin{itemize}
 \item
  $(\frac{p_{2}}{p_{1}})\neq 1$ and
 either $(\frac{-l}{p_{1}})=1$   {or} $l$$\equiv$$7$ (mod $8$);\\
 or
 \item
  $(\frac{p_{1}}{p_{2}})\neq 1$ and
 either $(\frac{-l}{p_{2}})=1$   {or} $l$$\equiv$$7$ (mod $8$).
 \end{itemize}
 \end{enumerate}
\end{proposition}
\begin{proof}
Simce $l$ is a prime, $l \equiv 3$ ({mod} $4$), applying Proposition \ref{twodotseven}, we have the following fields extension $\mathbb{Q}\left(i\sqrt{l}\right)\subset \mathbb{Q}\left(\xi_{l}\right).$\\
We have the following sequence of fields extensions
$$\mathbb{Q}\subset \mathbb{Q}\left(i\sqrt{l}\right)\subset \mathbb{Q}\left(\xi_{l}\right)= \mathbb{Q}\left(\xi^{l}_{l^{2}}\right)\subset \mathbb{Q}\left(\xi_{l^{2}}\right)\subset \mathbb{Q}\left(\xi_{l^{3}}\right)\subset \cdots
\subset \mathbb{Q}\left(\xi_{l^{k}}\right)$$
The degrees of the following field extensions are known: 
$$\left[\mathbb{Q}\left(i\sqrt{l}\right): \mathbb{Q}\right]=2, \left[\mathbb{Q}\left(\xi_{l}\right): \mathbb{Q}\right]=\varphi\left(l\right)=l-1,\ldots, \left[\mathbb{Q}\left(\xi_{l^{k}}\right): \mathbb{Q}\right]=\varphi\left(l^{k}\right)=l^{k-1}\cdot\left(l-1\right),$$
for any positive integer $k.$\\
Applying the transitivity theorem of finite extensions and taking into account the fact that $l \equiv 3$ ({mod} $4$), we obtain that the degree of the following field extensions are odd, namely:
$$\left[\mathbb{Q}\left(\xi_{l}\right):\mathbb{Q}\left(i\sqrt{l}\right)\right]=\frac{l-1}{2},
 \left[\mathbb{Q}\left(\xi^{2}_{l}\right): \mathbb{Q}\left(\xi_{l}\right)\right]=l,\ldots \left[\mathbb{Q}\left(\xi^{k}_{l}\right): \mathbb{Q}\left(\xi^{k-1}_{l}\right)\right]=l,$$
for any positive integer $k.$\\
Applying  Theorem \ref{theorem3.1} and Remark \ref{remark3.2}, we obtain that  the quaternion algebra $H_{\mathbb{Q}\left(\xi_{l^{k}}\right)}\left(p_{1}, p_{2}\right)$ is a division algebra if and only if one of the conditions 1., 2., 3. from the conclusion of Proposition \ref{fourdotone}
is fulfilled.

\end{proof}

\begin{corollary}
\label{fourdottwo}  Let $l$   {be an odd prime integer,} $l \equiv 3$ ({mod} $4$) and  let  $K=\mathbb{Q}\left(\xi_{l^{k}}\right)$ be the $l^{k}$th cyclotomic field, where $k$ is a positive integer.  Let $p_{1}, p_{2}$ be two distinct prime integers.

\begin{enumerate}[\rm1.]

 \item
Let the Kummer field $\mathbb{Q}\left(\xi_{l^{k}}, \sqrt[l^{k}]{\alpha}\right),$ where $\alpha$$\in$$K\backslash \{0\},$ $\alpha$ is different from the $l^{k}$th power of an integer in $K$. 
Then the quaternion algebra $H_{\mathbb{Q}\left(\xi_{l^{k}}, \sqrt[l^{k}]{\alpha}\right)}\left(p_{1}, p_{2}\right)$ is a division algebra if and only if 
one of the conditions from the conclusion of Proposition \ref{fourdotone}   holds;

\item
Let $L$ be any field extension of odd degree over $K.$ Then the quaternion algebra $H_{L}\left(p_{1}, p_{2}\right)$ is a division algebra if and only if 
one of the conditions from the conclusion of Proposition \ref{fourdotone}  holds.

 \end{enumerate}

\end{corollary}

\begin{proof}
Applying Proposition \ref{fourdotone} and Remark \ref{remark3.2}, the proof is immediate.
\end{proof}

Now we immediately obtain a result about quaternion algebras of the same type, over the $n$th cyclotomic field $\mathbb{Q}\left(\xi_{n}\right)$, when $n=2\cdot l^{k},$  with $l$  is a prime integer
$l\equiv 3$ (mod $4$).

\begin{corollary}
\label{fourdotthree}  Let $l$   {be an odd prime integer,} $l \equiv 3$ ({mod} $4$) and  let  $K=\mathbb{Q}\left(\xi_{n}\right)$ be the $n$th cyclotomic field, where $n=2\cdot l^{k}$ and  $k$ is a positive integer.  Let $p_{1}, p_{2}$ be two distinct prime integers. 
Then the quaternion algebra $H_{K}\left(p_{1}, p_{2}\right)$ is a division algebra if and only if  one of the conditions from the conclusion of Proposition \ref{fourdotone}  holds.

\end{corollary}

\begin{proof}
According to Proposition \ref{twodoteight} we have $\mathbb{Q}\left(\xi_{n}\right)= \mathbb{Q}\left(\xi_{l^{k}}\right).$  Applying Proposition \ref{fourdotone}  the proof is immediate.
\end{proof}

\begin{equation*}
\end{equation*}
\section{ Division quaternion algebras over a cyclotomic field  $\mathbb{Q}\left(\xi_{l}\right),$  when $l$  is a prime integer,
$l=2^{k}+1$ } 
\label{five}
\begin{equation*}
\end{equation*}

Let $l$  be a prime integer, $l\equiv 1$ (mod $4$) and let $p_{1},$  $p_{2}$ be two distinct prime integers. For to study when the quaternion  $H_{\mathbb{Q}\left(\xi_{l}\right)}\left(p_{1}, p_{2}\right)$ splits or it is a division algebra, we use Proposition  \ref{Proposition2.9} and the local-global principle.
Let $\Omega$ be the set of $p-$ adic completions $\mathbb{Q}_{p}$ of $\mathbb{Q},$ where $p$ is an arbitrary prime integer or $p=\infty$ (because $\mathbb{R}$$=\mathbb{Q}_{\infty}$). We have  $H_{K}\left(p_{1}, p_{2}\right)$ splits if and only if $H_{K_{p}}\left(p_{1}, p_{2}\right)$  splits over each of the completations $K_{p}$ (with respect the $p$-adic valuation) of $K=\mathbb{Q}\left(\xi_{l}\right).$ Since $K$ has not real embeddings, it results that we never need to consider $\mathbb{Q}_{\infty}$ for this problem on quaternion algebras.\\
Let $S$ be be the set of completions $\mathbb{Q}_{p}$$\in$$\Omega$ such that $H_{\mathbb{Q}}\left(p_{1}, p_{2}\right)$ is a division algebra (i.e a non-split algebra, according to Proposition \ref{twodotfive}) over  $\mathbb{Q}_{p}$. Applying Proposition \ref{Proposition2.9} and Corolarry \ref{cases 2}, it results that $S$ is a finite set and the cardinal of $S$ is an even number.\\
We denote by $\Omega^{'}$ the set of completions  $K_{\bf{p}}$ of the cyclotomic field $K=\mathbb{Q}\left(\xi_{l}\right).$ For each $\mathbb{Q}_{p}$$\in$$S$ and for each $K_{\bf{p}}$$\in$$\Omega^{'}$, where $\bf{p}$ $|$ $p,$ according to the local-global theorem we obtain that $H_{K}\left(p_{1}, p_{2}\right)$ splits 
if and only if $H_{K_{\bf{p}}}\left(p_{1}, p_{2}\right)$ splits. Therefore, in agreement with the above and applying Proposition \ref{twodotfive}, the quaternion algebra $H_{K}\left(p_{1}, p_{2}\right)$ is a division algebra
if and only if there are $\mathbb{Q}_{p}$$\in$$S$ and  $K_{\bf{p}}$$\in$$\Omega^{'}$, where\\
$\bf{p}$ $|$ $p$ such that  $H_{K_{\bf{p}}}\left(p_{1}, p_{2}\right)$ is a division algebra. Thus there are only finitely many conditions to check.
Applying Proposition \ref{Proposition2.15} and  Proposition \ref{twodotfive} we obtain that $H_{K_{\bf{p}}}\left(p_{1}, p_{2}\right)$ is a division algebra if and only if there are $\mathbb{Q}_{p}$$\in$$S$ and  $K_{\bf{p}}$$\in$$\Omega^{'},$ where\\
$\bf{p}$ $|$  $p$ such that  the degree $\left[K_{\textbf{p}}:\mathbb{Q}_{p}\right]$ is an odd number.  \ \ \ \ \ \ \ \ \ \ \ \ \ \ \ \ \ \ \ \ \ \ \ \ \ \ \ \ \ \ \ \ \ \ \ \ \ \ \ \       (5.1)\\

We obtain the following result.

\begin{proposition}
\label{fivedotone} Let $k$ be an integer, $k\geq 2$ such that $l=2^{k}+1$ is a prime number.
Let $\xi_{l}$ be a primitive root of order $l$ of the unity and  let  $\mathbb{Q}\left(\xi_{l}\right)$ be the $l$th cyclotomic field. Let $p_{1}, p_{2}$ be two distinct prime integers.
Then the quaternion algebra $H_{\mathbb{Q}\left(\xi_{l}\right)}\left(p_{1}, p_{2}\right)$ is a division algebra if and only if 
$p_{1}$ and $p_{2}$  are odd,  $p_{1}$$ \equiv 1$ (mod $l$), $p_{2}\neq l$ and $(\frac{p_{2}}{p_{1}})=-1$.
 \end{proposition}  

\begin{proof}
Since $l$ is prime  $l=2^{k}+1$, it immediately follows that $k=2^{n}$ (where $n$ is a positive integer), so $l$ is a Fermat prime.\\
Applying  (2.1) and (5.1),   $H_{\mathbb{Q}\left(\xi_{l}\right)}\left(p_{1}, p_{2}\right)$  is a division algebra if and only if there is a prime ideal $P$ of the ring $\mathbb{Z}\left[\xi_{l}\right],$
$P| p_{1}\mathbb{Z}\left[\xi_{l}\right]$ or  $P|p_{2}\mathbb{Z}\left[\xi_{l}\right]$ or $P|2\mathbb{Z}\left[\xi_{l}\right]$ such that $\left[K_{P}:\mathbb{Q}_{p}\right]$ is an odd number, where $p= p_{1}$ or $p= p_{2}$ or $p= 2.$\\
But, according to Proposition  \ref{*},  $\left[K_{P}:\mathbb{Q}_{p}\right]=e\cdot f$ where $f$ is the residual degree of $P$ in the fields extension $K\subset  L$.\\
If $p=l$ ( that is $p_{1}=l$ or $p_{2}= l$), applying Corollary \ref{2.5}, $p=l$ is ramified, namely  $p\mathbb{Z}\left[\xi_{l}\right]=l\mathbb{Z}\left[\xi_{l}\right]=\left(1- \xi_{l}\right)^{l-1} \mathbb{Z}\left[\xi_{l}\right]$  then $e=l-1,$  $g=1,$ $f=1$.
So $\left[K_{P}:\mathbb{Q}_{p}\right]=e\cdot f=l-1,$ so $\left[K_{P}:\mathbb{Q}_{l}\right]$ is not odd.\\
If $p=2,$ let  $f^{'}$ be the smallest positive integer such that $2^{f^{'}}_{1}\equiv 1$  (mod $l$), that is $f^{'}= ord_{\left(\mathbb{Z}^{*}_{l}, \cdot\right)}\left(\overline{2}\right) $ . According to a Corollary of Lagrange's Theorem, $f^{'}| \phi\left(l\right)=2^{2^{n}},$ where $n$ is a positive integer. It is clear that $f^{'}\neq 1$,
so   $f^{'}$  is even. Since  $\mathbb{Q}\subset \mathbb{Q}\left(\xi_{l}\right) $ is a Galois extension, according to Proposition \ref{*}, the ideal $2\mathbb{Z}\left[\xi_{l}\right]$
is written uniquely in the form
$$2\mathbb{Z}\left[\xi_{l}\right]=P^{e}_{1}\cdot P^{e}_{2}\cdot...\cdot P^{e}_{g},$$
where $g, e$ are positive integers and $P_{1}, P_{2},...., P_{g}$  are distinct prime ideals of the ring $\mathbb{Z}\left[\xi_{l}\right]$ and $e\cdot f \cdot g=\phi\left(l\right)=l-1,$ where  $f= \left[\mathbb{Z}\left[\xi_{l}\right]/ P_{i} : \mathbb{Z}/2 \mathbb{Z}\right]$ (for any $i=\overline{1,g}$  is the residual degree.
Applying Theorem \ref{theorem2.2}, we have $e=1$ and   $g=\frac{\varphi \left( l\right) }{f^{'}}.$ So $f=f^{'}$ and it is an even number. But, according to Proposition  \ref{*},  $\left[K_{P}:\mathbb{Q}_{2}\right]=e\cdot f.$ From what has been proved, it follows that $\left[K_{P}:\mathbb{Q}_{2}\right]= f$ is not odd.\\
Considering that S is a finite set whose cardinal is even and applying (5.1), we obtain that the quaternion algebras $H_{\mathbb{Q}\left(\xi_{l}\right)}\left(l, p_{2}\right)$ (with $p_{2}$ an odd number,  $p_{2}\neq l$), $H_{\mathbb{Q}\left(\xi_{l}\right)}\left(2, p_{l}\right)$ (with $p_{2}$ an odd number,  $p_{2}\neq l$),  $H_{\mathbb{Q}\left(\xi_{l}\right)}\left(2, l\right)$ are not division algebras. Taking into account these things, the definition of the set $S$ and
the fact that if  $H_{\mathbb{Q}\left(\xi_{5}\right)}\left(p_{1}, p_{2}\right)$ is a division algebra, it results $H_{\mathbb{Q}}\left(p_{1}, p_{2}\right)$ is a division algebra and appylying Lemma \ref{twodottwelve} and (2.1), it results we have one of the following situations: a) $p_{1}$ and $p_{2}$  are odd, $p_{1}$ or $p_{2} \equiv 1$ (mod $4$) and
$(\frac{p_{2}}{p_{1}})=-1$; b) $p_{1}$ and $p_{2}$  are odd, $p_{1}\equiv p_{2} \equiv 3$ (mod $4$) and
$(\frac{p_{2}}{p_{1}})=-1$. Reuning these two situations, we have: if $H_{\mathbb{Q}\left(\xi_{5}\right)}\left(p_{1}, p_{2}\right)$ if a division algebra, it follows that $p_{1}$ and $p_{2}$ are odd  $(\frac{p_{2}}{p_{1}})=-1$.\\

The last case that comes into the discussion  for to see if we obtain a division quaternion algebra  $H_{\mathbb{Q}\left(\xi_{l}\right)}\left(p_{1}, p_{2}\right)$ is:  when $p_{1}$ and $p_{2}$  are odd,  $p_{1}\neq l,$ $p_{2}\neq l,$\\
If  $p=p_{1}\neq l,$ $p_{2}\neq l,$   taking into account the above,  $H_{\mathbb{Q}\left(\xi_{l}\right)}\left(p_{1}, p_{2}\right)$ is a division algebra if and only if $\left[K_{P}:\mathbb{Q}_{p_{1}}\right]=e\cdot f$ is an odd number.\\
Since $p_{1}$ is an odd prime positive  integer and $\mathbb{Q}\subset \mathbb{Q}\left(\xi_{l}\right) $ is a Galois extension, according to Proposition \ref{*}, the ideal $p_{1}\mathbb{Z}\left[\xi_{l}\right]$
is written uniquely in the form
$$p_{1}\mathbb{Z}\left[\xi_{5}\right]=P^{e}_{11}\cdot P^{e}_{12}\cdot...\cdot P^{e}_{1g},$$
where $g, e$ are positive integers and $P_{11}, P_{12},...., P_{1g}$  are distinct prime ideals of the ring $\mathbb{Z}\left[\xi_{l}\right].$\\
Since $\mathbb{Q}\subset \mathbb{Q}\left(\xi_{l}\right) $ is a Galois extension, it results that $\left[K: \mathbb{Q}\right]=efg=l-1,$ where $f= \left[\mathbb{Z}\left[\xi_{l}\right]/ P_{1i} : \mathbb{Z}/p_{1} \mathbb{Z}\right]$ (for any $i=\overline{1,g}$) is the residual degree.\\
Applying Theorem \ref{theorem2.2} we have $e=1$ and $g=\frac{\varphi \left( l\right) }{f^{'}}=\frac{l-1}{f^{'}},$ where $f^{'}$ is the smallest positive integer such that $p^{f^{'}}_{1}\equiv 1$  (mod $l$).\\
 From all these, we obtain that is $\left[K_{P_{1}}:\mathbb{Q}_{p}\right]=f=f^{'}= ord_{\left(\mathbb{Z}^{*}_{l}, \cdot\right)}\left(\overline{p_{1}}\right) $ is and odd number. But according to a Corollary of Lagrange's Theorem we have
 $$ord_{\left(\mathbb{Z}^{*}_{l}, \cdot\right)}\left(\overline{p_{1}}\right) | \left(l-1\right) =2^{2^{n}},$$ 
where $n$ is a positive integer.\\
We obtain that $f^{'}=1,$ therefore $g=l-1,$ therefore $p_{1}$ is totally split in  $\mathbb{Z}\left[\xi_{l}\right],$ which is equivalent to $p_{1}$$ \equiv 1$ (mod $l$).\\
Combining these results, we obtain that the quaternion algebra  $H_{\mathbb{Q}\left(\xi_{l}\right)}\left(p_{1}, p_{2}\right)$ is a division algebra if and only if $p_{1}$ and $p_{2}$  are odd,  $p_{1}$$ \equiv 1$ (mod $l$), $p_{2}\neq 5$ and $(\frac{p_{2}}{p_{1}})=-1$.\\
The case when $p=p_{2}$ is similar.
\end{proof}

We remark that applying Proposition \ref{fivedotone} , for $l=5$ (that is, working in the $5$th cyclotomic field $\mathbb{Q}\left(\xi_{5}\right)$)  the conditions from  Proposition \ref{threedotnine} are not only sufficient, but necessary and sufficient. 

\begin{corollary}
\label{fivedottwo}
 Let $k$ be an integer, $k\geq 2$ such that $l=2^{k}+1$ is a prime number.
Let $\xi_{2l}$ be a primitive root of order $2l$ of the unity and  let  $\mathbb{Q}\left(\xi_{2l}\right)$ be the $2l$th cyclotomic field. Let $p_{1}, p_{2}$ be two distinct prime integers.
Then the quaternion algebra $H_{\mathbb{Q}\left(\xi_{2l}\right)}\left(p_{1}, p_{2}\right)$ is a division algebra if and only if 
$p_{1}$ and $p_{2}$  are odd,  $p_{1}$$ \equiv 1$ (mod $l$), $p_{2}\neq l$ and $(\frac{p_{2}}{p_{1}})=-1$.
   
\end{corollary}

\begin{proof} 

According to Proposition \ref{twodoteight} and Proposition \ref{fivedotone} the proof results immediately.

\end{proof}

We notice that applying Corollary \ref{fivedottwo}, for $l=5$ the conditions from  Proposition \ref{threedotten} are not only sufficient, but necessary and sufficient.

\begin{center}
\end{center}
\textbf{Acknowledgements.}
The author wants to thank Professor Victor Alexandru and Professor  Mirela \c Stef\u anescu, for the useful discussions on
this topic.

\begin{center}
\end{center}

\begin{equation*}
\end{equation*}

\end{document}